\documentstyle[amscd]{amsart}
\numberwithin{equation}{section}
\theoremstyle{plain}
\newtheorem{thm}{Theorem}[section]

\newtheorem{lem}[thm]{Lemma}
\newtheorem{prop}[thm]{Proposition}
\begin{document}
\title{K-theory for the simple $C^*$-algebra \\
of the Fibonacchi Dyck system}
\author{Kengo Matsumoto}
\address{ 
Department of Mathematical Sciences, 
Yokohama City  University,
22-2 Seto, Kanazawa-ku, Yokohama, 236-0027 Japan}
\email{kengo@@yokohama-cu.ac.jp}
\maketitle
\begin{abstract}
Let $F$
be the Fibonacci matrix
$
\bigl[\begin{smallmatrix}
1 & 1 \\
1 & 0 \\
\end{smallmatrix}\bigr] 
$.
The Fibonacci Dyck shift is a subshsystem of the Dyck shift $D_2$
constrained by the matrix $F$.
Let ${{\frak L}^{Ch(D_F)}}$ 
be a $\lambda$-graph system presenting the subshift
$D_F$, that is called the Cantor horizon $\lambda$-graph system
for $D_F$.
We will study the
 $C^*$-algebra 
${\cal O}_{{\frak L}^{Ch(D_F)}}$ 
associated with 
$
{{\frak L}^{Ch(D_F)}}
$.
It is simple purely infinite and generated by 
four partial isometries with some operator relations.
We will compute the K-theory of the $C^*$-algebra. 
As a result, the $C^*$-algebra is simple purely infinite and not semiprojective. Hence it is not stably isomorphic to any  
Cuntz-Krieger algebra.
\end{abstract}

\def\Zp{{ {\Bbb Z}_+ }}
\def\U{{ {\cal U} }}
\def\S{{ {\cal S} }}
\def\M{{ {\cal M} }}
\def\OL{{ {\cal O}_\Lambda }}
\def\A{{ {\cal A} }}
\def\Ext{{{\operatorname{Ext}}}}
\def\Im{{{\operatorname{Im}}}}
\def\Hom{{{\operatorname{Hom}}}}
\def\Ker{{{\operatorname{Ker}}}}
\def\dim{{{\operatorname{dim}}}}
\def\id{{{\operatorname{id}}}}
\def\OLF{{{\cal O}_{{\frak L}^{Ch(D_F)}}}}
\def\OLN{{{\cal O}_{{\frak L}^{Ch(D_N)}}}}
\def\OLA{{{\cal O}_{{\frak L}^{Ch(D_A)}}}}
\def\LCHDA{{{{\frak L}^{Ch(D_A)}}}}
\def\LCHDN{{{{\frak L}^{Ch(D_N)}}}}
\def\LCHDF{{{{\frak L}^{Ch(D_F)}}}}
\def\LCHLA{{{{\frak L}^{Ch(\Lambda_A)}}}}
\def\LWA{{{{\frak L}^{W(\Lambda_A)}}}}





\section{Introduction}

Let $\Sigma$ be a finite set with its discrete topology, 
that is called an alphabet.
Each element of $\Sigma$ is called a symbol.
Let $\Sigma^{\Bbb Z}$ 
be the infinite product space 
$\prod_{i=-{\infty}}^{\infty}\Sigma_{i},$ 
where 
$\Sigma_{i} = \Sigma$,
 endowed with the product topology.
 The transformation $\sigma$ on $\Sigma^{\Bbb Z}$ 
given by 
$\sigma({(x_i)}_{i \in \Bbb Z}) = (x_{i+1})_{i \in \Bbb Z}$ 
is called the full shift over $\Sigma$.
 Let $\Lambda$ be a closed subset of 
 $\Sigma^{\Bbb Z}$ such that $\sigma(\Lambda) = \Lambda$. 
 The topological dynamical system 
  $(\Lambda, \sigma\vert_{\Lambda})$
is called a subshift or a symbolic dynamical system.
 It is written as $\Lambda$ for brevity.
 
In \cite{Ma}, 
the  author  
has introduced a notion of 
$\lambda$-graph system 
as  a presentation of subshifts.
A $\lambda$-graph system 
$ {\frak L} = (V,E,\lambda,\iota)$
consists of a vertex set 
$V = V_0\cup V_1\cup V_2\cup\cdots$, an edge set 
$E = E_{0,1}\cup E_{1,2}\cup E_{2,3}\cup\cdots$, 
a labeling map
$\lambda: E \rightarrow \Sigma$
and a surjective map
$\iota_{l,l+1}: V_{l+1} \rightarrow V_l$ for each
$l\in \Zp$, 
where $\Zp$ denotes the set of all 
nonnegative integers.
An edge $e \in E_{l,l+1}$ has its source vertex $s(e)$ in $V_l,$
its terminal vertex $t(e)$ in $V_{l+1}$ 
and its label $\lambda(e)$ in $\Sigma$ (\cite{Ma}).

The theory of symbolic dynamical system has  a close relationship to automata theory and language theory.
In the theory of language, there is a class of universal languages due to W. Dyck.
The symbolic dynamics generated by the languages are called the Dyck shifts 
$D_N$ (cf. \cite{ChS}, \cite{Kr},\cite{Kr2},\cite{Kr3}). 
  Its alphabet consists of the $2N$ brackets: 
$(_1, \dots,(_N, )_1, \dots,)_N$.
 The forbidden words consist of words that do not obey the standard bracket rules. 
In \cite{KM},  
a $\lambda$-graph system ${\frak L}^{Ch(D_N)}$ 
that presents the subshift $D_N$
has been introduced.
The $\lambda$-graph system is called the Cantor horizon  
$\lambda$-graph system 
for the Dyck shift 
$D_N$.
The K-groups for ${\frak L}^{Ch(D_N)}$,
that are invariant under topological conjugacy of the subshift $D_N$,
have been computed (\cite{KM}). 

In \cite{Ma6} (cf. \cite{KM}), 
the $C^*$-algebra $\OLN$ associated with the 
Cantor horizon $\lambda$-graph system $\LCHDN$ has been studied.
In the paper, it has been proved that the algebra 
$\OLN$ is simple and purely infinite and generated by  
$N$ partial isometries and $N$ isometries satisfying some operator relations.
Its K-groups are
$$
K_0({\cal O}_{{\frak L}^{Ch(D_N)}}) \cong 
{\Bbb Z}/N{\Bbb Z} \oplus C({\frak K},{\Bbb Z}),
\qquad
K_1({\cal O}_{{\frak L}^{Ch(D_N)}}) \cong 0
$$
where $C(\frak K,{\Bbb Z})$ denotes the abelian group of all integer valued continuous functions on a Cantor discontinuum $\frak K$
(cf. \cite{KM}).

Let $u_1,\dots,u_N$ be 
the canonical generating isometries 
of the Cuntz algebra ${\cal O}_N$ 
that satisfy the relations:
$
\sum_{j=1}^{N} u_ju_j^* = 1, \,     u_i^*u_i = 1 
$
for 
$ i=1,\dots N. 
$
Then the bracket rule of the symbols $(_1, \dots,(_N,\,  )_1,\dots, )_N$
of the Dyck shift $D_N$ 
may be interpreted as the relations
$
u_i^* u_i = 1, \,
u_i^* u_j = 0 
$ 
for
$ i \ne j
$
of the partial isometries  
$u_1^*, \dots, u_N^*, u_1, \dots, u_N$
in the $C^*$-algebra ${\cal O}_N$ (cf. (2.1)).

In \cite{Ma7}, 
we have considered
a generalization of Dyck shifts $D_N$ 
by using the canonical generators of the Cuntz-Krieger algebras 
${\cal O}_A$
 for  
$N \times N$ matrices $A$ with entries in $\{ 0,1 \}$.
The generalized Dyck shift is denoted by $D_A$ 
and called the topological Markov Dyck  shift for $A$ (cf. \cite{HIK}, \cite{KM2}).
Let $\alpha_1, \dots,\alpha_N, \, \beta_1, \dots, \beta_N$
be the alphabet of $D_A$.
They correspond to the brackets 
$(_1, \dots,(_N,\,  )_1,\dots, )_N$
respectively.
Let $t_1, \dots, t_N$ be 
the canonical generating partial isometries 
of the Cuntz-Krieger  algebra ${\cal O}_A$ 
that satisfy the relations:
$
\sum_{j=1}^{N} t_jt_j^* = 1, \,     
t_i^*t_i = \sum_{j=1}^N A(i,j)t_jt_j^*
$
for 
$ i=1,\dots, N. 
$
Consider the correspondence
$
\varphi(\alpha_i) = t_i^*,
\varphi(\beta_i) = t_i, i=1,\dots,N.
$
Then a word $w$ of $\{\alpha_1, \dots, \alpha_N, \beta_1,\dots,\beta_N \}$
is defined to be admissible for the subshift $D_A$ 
precisely  if
the correspnoding element to $w$ through $\varphi$ in ${\cal O}_A$
is not zero.
Hence we may recognize  $D_A$ to be  the subshift defined by the canonical generators of the Cuntz-Krieger algebra ${\cal O}_A$.  
The subshifts  $D_A$ are not sofic in general and reduced to the Dyck shifts if all entries of $A$ are $1$.

The Cantor horizon 
$\lambda$-graph system $\LCHDA$
for the topological Markov Dyck shift $D_A$
has been also studied in \cite{Ma7}.
It has been proved  to be 
$\lambda$-irreducible with $\lambda$-condition (I) in the sense of \cite{Ma5}
if the matrix is irreducible with condition (I) in the sense of Cuntz-Krieger \cite{CK}.
Hence the  associated $C^*$-algebra $\OLA$ is simple and purely infinite.
It is the unique $C^*$-algebra generated by $2N$ partial isometries 
subject to some operator relations.

 In this paper we  study the $C^*$-algebra
 $\OLF$ for the Fibonacci matrix 
$ F = 
\bigl[\begin{smallmatrix}
1 & 1 \\
1 & 0 \\
\end{smallmatrix}\bigr]
$.
It is  the smallest matrix in the irreducible square matices with condition (I) such that the associated topological Markov shift $\Lambda_F$
is not conjugate to any full shift.
The topological entropy of 
$\Lambda_F$  is
$\log \frac{1 + \sqrt{5}}{2}$ 
the logarithm of the Perron eigenvalue of $F$.     
We call the subshift $D_F$ the Fibonacci Dyck shift.
As the matrix is irreducible with condition (I),
the associated $C^*$-algebra
$\OLF$ is simple and purely infinite.
We will compute the K-groups 
$K_i(\OLF), i=0,1$ of the algebra
so that we have
\begin{thm}
The $C^*$-algebra $\OLF$
associated with the $\lambda$-graph system
$\LCHDF$ is  unital, separable, nuclear, simple and purely infinite.
It is the unique $C^*$-algebra generated by one isometry $T_1$
and three  partial isometries 
$S_1, S_2, T_2$  
subject to the following operator relations:
{\allowdisplaybreaks
\begin{align}
\sum_{j=1}^{2} & ( S_jS_j^* +   T_jT_j^* ) 
  =  \sum_{j=1}^{2}  S_j^*S_j      = 1,  \qquad 
T_2^*T_2   =   S_1^*S_1,  \\
E_{\mu_1\cdots \mu_k} & = \sum_{j=1}^2 F(j,\mu_1) S_jS_j^*
E_{\mu_1\cdots \mu_k}
S_jS_j^* + T_{\mu_1}E_{\mu_2\cdots \mu_k}T_{\mu_1}^*, \qquad k >1    
\end{align}}
where 
$E_{\mu_1\cdots \mu_k}
= S_{\mu_1}^*\cdots S_{\mu_k}^*S_{\mu_k}\cdots S_{\mu_1}$,
$(\mu_1,\cdots,\mu_k)\in \Lambda_F^*,$ 
and $\Lambda_F^*$ is the set of admissible words of the topological Markov shift $\Lambda_F$ defined by the matrix $F$.
The K-groups are
$$
K_0(\OLF) \cong {\Bbb Z} \oplus  C({\frak K},{\Bbb Z})^{\infty},
\qquad
K_1(\OLF) \cong 0.
$$
\end{thm}
This paper is a continuation of \cite{Ma7}.

\section{The subshift $D_A$ and the $\lambda$-graph system $\LCHDA$}
We will briefly review the topological Markov Dyck shift $D_A$ and 
its Cantor horizon $\lambda$-graph system $\LCHDA$.

Consider a pair of $N$ symbols 
where
$
\Sigma^- = \{ \alpha_1,\cdots,\alpha_N \},
\Sigma^+ = \{ \beta_1,\cdots,\beta_N \}.
$
We set  
$\Sigma = \Sigma^- \cup \Sigma^+$.
Let $A =[A(i,j)]_{i,j=1,\dots,N}$
be an $N\times N$ matrix with entries in $\{0,1\}$.
Throughout this paper,
 $A$ is assumed to have no zero rows or columns.
Consider the Cuntz-Krieger algebra ${\cal O}_A$ for the matrix $A$
that is the universal $C^*$-algebra generated by 
$N$ partial isometries $t_1,\dots,t_N$ subject to the following relations:
$$
\sum_{j=1}^N t_j t_j^* = 1, 
\qquad
t_i^* t_i = \sum_{j=1}^N A(i,j) t_jt_j^* \quad \text{ for } i = 1,\dots,N
$$
(\cite{CK}).
Define a correspondence 
$\varphi_A :\Sigma \longrightarrow \{t_1^*,\dots, t_N^*, t_1,\dots, t_N \}$
by setting
$$ 
\varphi_A(\alpha_i) = t_i^*,\qquad 
\varphi_A(\beta_i) = t_i  \quad \text{ for } i=1,\dots,N.
$$
We denote by $\Sigma^*$ the set of all words 
$\gamma_1\cdots \gamma_n$ of elements of $\Sigma$.
Define the set
$$
{\frak F}_A = \{ \gamma_1\cdots \gamma_n \in \Sigma^* \mid
\varphi_A(\gamma_1)\cdots \varphi_A( \gamma_n) = 0 \}.
$$
Let $D_A$ be the subshift over $\Sigma$ whose forbidden words are 
${\frak F}_A.$
The subshift is called the topological Markov Dyck shift defined by $A$
(cf. \cite{HIK}, \cite{KM2}).
If all entries of $A$  are $1$, 
the subshift $D_A$ 
becomes the Dyck shift $D_N$
with $2N$ bracket
 (cf. \cite{Kr2},\cite{Kr3}, \cite{KM}, \cite{Ma6},\cite{Ma7}).
We note the fact that 
 $\alpha_i \beta_j\in {\frak F}_A$  if $i\ne j$,
 and
 $\alpha_{i_n}\cdots \alpha_{i_1} \in {\frak F}_A$  
if and only if 
$\beta_{i_1}\cdots \beta_{i_n} \in {\frak F}_A$.
Consider the following  subsystem of $D_A$
$$
D_A^+  = \{ {(\gamma_i)}_{i \in \Bbb Z} \in D_A \mid
\gamma_i \in \Sigma^+ \text{ for all } i \in \Bbb Z \}.
$$
The subshift  
$D_A^+$ is identified with the topological Markov shift 
$$
\Lambda_A = \{ {(x_i)}_{i \in \Bbb Z}\in \{ 1,\dots,N \}^{\Bbb Z}
 \mid A(x_i,x_{i+1}) = 1, i \in \Bbb Z \}
$$ 
defined by the matrix $A$.
Hence the subshift $D_A$ is recognized to contain the topological Markov shift 
$\Lambda_A$.

We denote by 
$B_l(D_A)$
and 
$B_l(\Lambda_A)$
 the set of admissible words of length 
$l$ of $D_A$
and that of 
$\Lambda_A$ respectively.
Let $m(l)$ be the cardinal number of $B_l(\Lambda_A)$.
We use lexcographic order from the left on the words of $B_l(\Lambda_A)$,
so that we may assign to a word $\mu_1\cdots \mu_l\in B_l(\Lambda_A)$ 
the number $N(\mu_1\cdots \mu_l)$ from $1$ to $m(l)$.
For example, 
if 
$A = F=
\bigl[\begin{smallmatrix} 
1 & 1 \\
1 & 0 \\
\end{smallmatrix}
\bigr],
$
then 
\begin{align*}
B_1(\Lambda_F)&  = \{ 1,2\},\qquad N(1) =1, \, N(2) = 2, \\
B_2(\Lambda_F)& = \{ 11, 12, 21 \},\qquad N(11) =1,\,  N(12) = 2,\, N(21) = 3,
\end{align*}
and so on.
Hence the set $B_l(\Lambda_A)$  bijectively corresponds 
to the set of natural numbers less than or equal to $m(l)$.  
Let us now describe 
the Cantor horizon $\lambda$-graph system $\LCHDA$ of $D_A$.
The vertices $V_l$ at level $l$ for $l \in \Zp$
are given by the admissible words of length $l$
consisting of the symbols of $\Sigma^+$.
We regard $V_0$ as a one point set of the empty word $\{ \emptyset \}$. 
Since $V_l$ is identified with $B_l(\Lambda_A)$,
we may write $V_l$ as
$$
V_l = \{ v^l_{N(\mu_1 \cdots \mu_l)}
 \mid \mu_1\cdots\mu_l\in B_l(\Lambda_A) \}.
$$
The mapping $\iota ( = \iota_{l,l+1}) :V_{l+1}\rightarrow V_l$ 
is defined by deleting the rightmost symbol of a corresponding word such as 
$$ 
\iota( v^{l+1}_{N(\mu_1\cdots \mu_{l+1})}) 
= v^l_{N(\mu_1 \cdots \mu_l)}\quad \text{ for }\quad
v^{l+1}_{N(\mu_1 \cdots \mu_{l+1})} \in V_{l+1}.
$$
We define an edge labeled $\alpha_j$ from
$v^l_{N(\mu_1 \cdots \mu_l)}\in V_l$ 
to
$v^{l+1}_{N(\mu_0\mu_1 \cdots \mu_l)}\in V_{l+1}$
precisely if
$\mu_0 = j,$
and 
 an edge labeled $\beta_j$ from
$v^l_{N(j\mu_1 \cdots \mu_{l-1})}
\in V_l$ 
to
$v^{l+1}_{N(\mu_1 \cdots \mu_{l+1})}\in V_{l+1}.$
For $l=0$, 
we define an edge labeled $\alpha_j$
form $v^0_1$ to $v^1_{N(j)}$, 
and
an edge labeled $\beta_j$
form $v^0_1$ to $v^1_{N(i)}$ if $A(j,i) =1$. 
We denote by $E_{l,l+1}$ the set of edges from $V_l$ to $V_{l+1}$.
Set $E = \cup_{l=0}^{\infty} E_{l,l+1}$.
It is easy to see that the resulting labeled Bratteli diagram with 
$\iota$-map
becomes a $\lambda$-graph system over $\Sigma$, that is 
called the Cantor horizon $\Lambda$-graph system and is denoted by
$\LCHDA$.

A $\lambda$-graph system $\frak L$ is said to present a subshift $\Lambda$
if the set of all admissible words of $\Lambda$ coincides with the set of all 
finite labeled sequences appearing in concatenating edges of $\frak L$. 
In \cite{Ma7}, the following propositions have been proved.
\begin{prop} \hspace{5cm}
\begin{enumerate}
\renewcommand{\labelenumi}{(\roman{enumi})}
\item
If $A$ satisfies condition (I) in the sense of Cuntz-Krieger \cite{CK},
the subshift $D_A$ is not sofic.
\item
 The $\lambda$-graph system $\LCHDA$ presents the subshift $D_A$.
\item
If $A$ is an irreducible matrix with condition (I), then
the $\lambda$-graph system $\LCHDA$ is $\lambda$-irreducible with $\lambda$-condition (I) in the sense of \cite{Ma5}.
\end{enumerate}
\end{prop}

\begin{prop}
 The $C^*$-algebra $\OLA$
associated with the $\lambda$-graph system
$\LCHDA$ is unital, separable,  nuclear, simple and purely infinite.
It is the unique $C^*$-algebra generated by $2N$ partial isometries 
$S_i, T_i, i=1,\dots,N$  
subject to the following operator relations:
\begin{gather*}
\sum_{j=1}^{N}  ( S_jS_j^* +   T_jT_j^* ) = \sum_{j=1}^{N}  S_j^*S_j = 1,  \\
 T_i^*T_i  =  \sum_{j=1}^N A(i,j) S_j^*S_j, \qquad i=1,2,\dots, N, \\
E_{\mu_1\cdots \mu_k}  = \sum_{j=1}^N A(j,\mu_1) S_jS_j^*
E_{\mu_1\cdots \mu_k}
S_jS_j^* + T_{\mu_1}E_{\mu_2\cdots \mu_k}T_{\mu_1}^*, \qquad k >1    
\end{gather*}
where 
$E_{\mu_1\cdots \mu_k}
= S_{\mu_1}^*\cdots S_{\mu_k}^*S_{\mu_k}\cdots S_{\mu_1}$,
$(\mu_1,\cdots,\mu_k)\in \Lambda_A^*$ 
 the set of admissible words of the topological Markov shift $\Lambda_A$ defined by the matrix $A$.
\end{prop}


\section{K-theory for  $\OLF$}

We will prove Theorem 1.1.
The operator relations (1.1) and (1.2) are direct from the operator relations
in Proposition 2.2.
By Proposition 2.2, 
it remains to prove the K-group formulae.
This section is devoted to computing the K-groups
$K_i(\OLF), i=0,1$
 for the $C^*$-algebra $\OLF$.
 The symbols 
 $\alpha_1,\, \alpha_2,\, \beta_1,\, \beta_2$
 of the subshift $D_F$ correspond to the brackets
 $(_1,\, (_2,\, )_1,\, )_2$ respectively.
Let $V_l, l\in \Bbb Z$ be the vertex set of 
the $\lambda$-graph system $\LCHDF$.
They are identified with the admissible words consisting of 
the symbols $\beta_1, \beta_2$ in $\Sigma^+$.
Since the word $\beta_2\beta_2 $ is forbidden, the following is a list of the vertex sets $V_l$ for $l=0,1,2,3,4,\dots $:
\begin{align*}
V_0:  \medspace &  * \\
V_1:  \medspace &  (\beta_1),\, (\beta_2),\\
V_2:  \medspace &  (\beta_1\beta_1),\, (\beta_1\beta_2),\, (\beta_2\beta_1),\\
V_3:  \medspace &  (\beta_1\beta_1\beta_1),\, (\beta_1 \beta_1\beta_2),\, (\beta_1\beta_2\beta_1),\, (\beta_2 \beta_1\beta_1),\, (\beta_2\beta_1\beta_2),\\
V_4:  \medspace &  (\beta_1\beta_1\beta_1 \beta_1),\,  (\beta_1\beta_1\beta_1\beta_2),\, (\beta_1\beta_1 \beta_2\beta_1),\, (\beta_1\beta_2\beta_1\beta_1),\,
(\beta_1\beta_2\beta_1\beta_2),\,\\
& (\beta_2\beta_1\beta_1 \beta_1),\, (\beta_2\beta_1\beta_1\beta_2),\, (\beta_2\beta_1 \beta_2\beta_1),\\
& \cdots 
\end{align*}

Let $f_l$ be the $l$-th Fibonacci number for $l \in \Bbb N$.
They are inductively defined by  
$$
f_1 = f_2 = 1, \qquad  
f_{l+2} = f_{l+1} + f_l \quad \text{for }l \in \Bbb N.
$$
By the structure of the $\lambda$-graph system
$\LCHDF$,
the number $m(l)$ of the vertex set $V_l$ is $f_{l+2}$.
 We denote by 
 $(\M_{l,l+1}, I_{l,l+1})_{l\in \Zp}$
the symbolic matrix system 
of the Cantor horizon $\lambda$-graph system
$\LCHDF$.
We write the vertex set 
$V_l$ as $\{ v_1^l,\dots,v_{m(l)}^l \}.$
Both the matrices $\M_{l,l+1}$ and $I_{l,l+1}$ 
are the $m(l) \times m(l+1)$ matrices
for each $l \in \Zp.$
For $i=1,\dots,m(l), \, j=1,\dots,m(l+1)$,
the component $\M_{l,l+1}(i,j)$ denotes the formal sum of labels of edges starting at the vertex $v_i^l$ and terminating at the vertex $v_j^{l+1},$
and the component $I_{l,l+1}(i,j) $ denotes  $1$ if $\iota(v_j^{l+1} = v_i^l$,
otherwise $0$. 
They satisfy the  relations 
$
I_{l,l+1} \M_{l+1,l+2} = \M_{l,l+1}I_{l+1,l+2}
$ 
for 
$
l \in \Zp
$
as symbolic matrices.
The orderings of the rows and columns of the matrices are arranged lexcographically  on indices ${i_1}\cdots {i_n}$ of the words 
$\beta_{i_1}\cdots \beta_{i_n}$ from the left. 
Let us denote by 
$0_{p,q}$ the $m(p) \times m(q)$ matrix all of whose entries are  $0'$s. 
\begin{lem}
The  $m(l)\times m(l+1)$ matrix $I_{l,l+1}$ is  given by :
$$
I_{0,1} = 
\begin{bmatrix}
1&  1\\
\end{bmatrix},
\qquad
I_{1,2} 
=
\begin{bmatrix}
  1 & 1 & 0    \\
  0 & 0 & 1  
\end{bmatrix}, 
\qquad
I_{l+2,l+3} 
=
\begin{bmatrix} 
I_{l+1,l+2} & 0_{l+1,l+1}\\
0_{l,l+2}   & I_{l,l+1}
\end{bmatrix},
\quad l \in \Zp.
$$
\end{lem}
In what follows,
blanks at components of matrices denote 0's.
 For $l\in \Zp$ and 
$a \in \{\alpha_1,\alpha_2,\beta_1,\beta_2 \}$,
let $I_l(a)$ be the $m(l) \times m(l)$
diagonal matrix with diagonal entries $a,$
and $\S_l(a)$ the $m(l-1)\times m(l+1)$ 
matrix 
defined by
$$
\S_0(a) = 
\begin{bmatrix}
a & a \\
\end{bmatrix},
\quad
\S_1(a) =
\begin{bmatrix}
a & a & a \\
\end{bmatrix},
\quad
\S_{l+2}(a) =
\begin{bmatrix}
\S_{l+1}(a) & 0_{l,l+1} \\
0_{l-1,l+2}& \S_{l}(a) 
\end{bmatrix}
$$
where
$m(-1)$ denotes $1$.
For $l =2,3,4$  and
$a \in \{\alpha_1,\alpha_2,\beta_1,\beta_2 \}$,
one sees that
{\allowdisplaybreaks
\begin{align*}
{\S}_{2}(a) 
 &{= 
	\begin{bmatrix}
	a &a & a &   &        \\
  	  &  &   & a & a      
	\end{bmatrix}
 		\quad : 2 \times 5 \text{ matrix},}\\
{\S}_{3}(a) 
 &{ = 
	\begin{bmatrix}
	a& a& a&   &  &  &  &   \\
 	 &  &  & a &a &  &  &   \\
 	 &  &  &   &  & a&a &a  
	\end{bmatrix}
		\quad : 3 \times 8\text{ matrix},}
\end{align*}}
and
$$
\setcounter{MaxMatrixCols}{25}
{\S}_{4}(a) 
  = \begin{bmatrix}
	 a &a &a &  &  &  &  &  &  &  &  &  &  \\
	   &  &  &a &a &  &  &  &  &  &  &  &  \\
	   &  &  &  &  &a &a &a &  &  &  &  &  \\
	   &  &  &  &  &  &  &  &a &a &a &  &  \\
	   &  &  &  &  &  &  &  &  &  &  &a &a 
	\end{bmatrix}
		\quad : 5 \times 13\text{ matrix}.
$$
\begin{lem}
The $m(l) \times m(l+1)$ matrix $\M_{l,l+1}$ 
is given by :
\begin{align*}
\M_{0,1} & = 
\left[
\alpha_1 + \beta_1 + \beta_2,\,\,
\alpha_2 + \beta_1 
\right],\\
\M_{l,l+1} & = 
\left[\begin{array}{c|c}
\multicolumn{2}{c}{{\S}_l(\beta_1)}\\
\hline
{\S}_{l-1}(\beta_2)&  0_{l-2, l-1}
\end{array}\right]
 + 
\left[
\begin{array}{c|c}
                  &  I_{l-1}(\alpha_2) \\
{I}_l({\alpha_1}) & \hrulefill                  \\  
                  &  0_{l-2,l-1}
\end{array}\right]. 
\end{align*}
\end{lem}
\begin{pf} In the right hand side in the second equation above,
the first summand 
describes the transitions that arise when 
a vertex accepts a symbol in $\Sigma^{+}.$ 
The second summand describes the transitions that arise  when a vertex accepts a symbol in $\Sigma^-$.
\end{pf}
We present the above matrices for $ l =  1, 2, 3, 4$: 
{\small
\begin{gather*}
I_{1,2} 
=
\left[\begin{smallmatrix}
  1 & 1 &     \\
    &   & 1  
\end{smallmatrix}
\right],
\qquad
\M_{1,2}
=
\left[\begin{smallmatrix}
\alpha_1 + \beta_1 & \beta_1           & \alpha_2 + \beta_1  \\
\beta_2            &\alpha_1 + \beta_2 & 
\end{smallmatrix}
\right], \\
I_{2,3} 
= 
\left[\begin{smallmatrix}
 1&1 &  &  &   \\
  &  &1 &  &   \\
  &  &  &1 &1  
 \end{smallmatrix}\right], 
\qquad
\M_{2,3}=
\left[\begin{smallmatrix}
\alpha_1+\beta_1&\beta_1 & \beta_1        & \alpha_2 &                  \\
                &\alpha_1&                & \beta_1  &\alpha_2 + \beta_1\\
\beta_2         &\beta_2 &\alpha_1+\beta_2&          &                  \\     \ \end{smallmatrix}\right],\\
I_{3,4} 
=\left[\begin{smallmatrix}
1&1& & & & & &  \\
 & &1& & & & &  \\
 & & &1&1& & &  \\
 & & & & &1&1&  \\
 & & & & & & &1 \\
 \end{smallmatrix}\right],
\qquad
\M_{3,4}
=\left[\begin{smallmatrix}
\alpha_1
+\beta_1&\beta_1 &\beta_1 &       &       &\alpha_2&        &        \\                &\alpha_1&        &\beta_1&\beta_1&        &\alpha_2&        \\                &        &\alpha_1&       &       &\beta_1 &\beta_1 &\alpha_2 
                                                              + \beta_1 \\      \beta_2&\beta_2 & \beta_2&\alpha_1&       &        &        &          \\             &        &        &\beta_2 &\alpha_1
                                  + \beta_2&        &        &      \\
\end{smallmatrix}\right],\\
\setcounter{MaxMatrixCols}{25}
I_{4,5} 
=\left[\begin{smallmatrix}
1&1& & & & & & & & & & &  \\
 & &1& & & & & & & & & &  \\
 & & &1&1& & & & & & & &  \\
 & & & & &1&1& & & & & &  \\
 & & & & & & &1& & & & &  \\
 & & & & & & & &1&1& & &  \\
 & & & & & & & & & &1& & \\
 & & & & & & & & & & &1&1
\end{smallmatrix}\right],\\
\setcounter{MaxMatrixCols}{25}
\M_{4,5}
=\left[\begin{smallmatrix}
\alpha_1+
\beta_1 &\beta_1 & \beta_1&        &        &        &        &        &
\alpha_2&        &        &        &         \\  
        &\alpha_1&        &\beta_1 &\beta_1 &        &        &        &
        &\alpha_2&        &        &         \\
        &        &\alpha_1&        &        &\beta_1 &\beta_1 & \beta_1&
        &        &\alpha_2&        &         \\ 
        &        &        &\alpha_1&        &        &        &        &       \beta_1 &\beta_1 &\beta_1 &\alpha_2&         \\                                        &        &        &        &\alpha_1&        &        &        &               &        &        &\beta_1 &\beta_1 +\alpha_2    \\                     \beta_2&\beta_2 & \beta_2&        &        &\alpha_1&        &        &               &        &        &        &         \\           
        &        &        &\beta_2 &\beta_2 &        &\alpha_1&        &               &        &        &        &         \\            
        &        &        &        &        &\beta_2 &\beta_2 &\beta_2+\alpha_1&       &        &        &        &         \\
\end{smallmatrix}\right].
\end{gather*}}

Let $(M_{l,l+1}, I_{l,l+1})_{l\in \Zp}$ 
be the  nonnegative  matrix system for $(\M,I)$.
The matrix $M_{l,l+1}$ for each $l \in \Zp$
 is obtained from ${\M}_{l,l+1}$ by setting all the symbols of 
 ${\M}_{l,l+1}$ equal to $1$. 
That is, the $(i,j)$-component $M_{l,l+1}(i,j)$ of the matrix $M_{l,l+1}$ 
denotes the number of the symbols in $\Sigma$ that appear in 
$\M_{l,l+1}(i,j).$
The groups
$K_0(\OLF), K_1(\OLF)$
are realized as the K-groups $K_0(M,I)$ and $K_1(M,I)$
for the nonnegative matrix system $(M,I)$ respectively (cf. \cite{Ma2}).
They
are calculated by the following formulae. 
\begin{lem}[(\cite{Ma2}, cf. \cite{Ma})] 
\hspace{4cm}

\begin{enumerate}
\renewcommand{\labelenumi}{(\roman{enumi})}%
\item
$
K_0(\OLF ) 
=  
\underset{l }{\varinjlim} 
\{ {\Bbb Z}^{m(l+1)} / (M_{l,l+1}^t - I_{l,l+1}^t){{\Bbb Z}^{m(l)}}, 
   \bar{I}^t_{l,l+1} \},
$
where the inductive limit is taken along the natural induced homomorphisms
$
   \bar{I}^t_{l,l+1}, l \in \Zp
$
by the 
matrices $I^t_{l,l+1}$. 
\item
$
K_1(\OLF)
=\underset{l}{\varinjlim}\{ \Ker (M_{l,l+1}^t - I_{l,l+1}^t)
             \text{ in }  {\Bbb Z}^{m(l)}, I^t_{l,l+1} \}, 
$
where the inductive limit is taken along the homomorphisms of the restrictions of $I^t_{l,l+1}$ to  $\Ker (M_{l,l+1}^t - I_{l,l+1}^t).$
\end{enumerate}
\end{lem}
By the formulae of $\M_{l,l+1}$
in Lemma 3.2,
the matrices $M_{l,l+1}^t - I_{l,l+1}^t $
for $l=1,2,3,4$
are presented as in the following way:
{\allowdisplaybreaks
\begin{gather*}
{M_{1,2}^t - I_{1,2}^t=
\left[\begin{smallmatrix}
\begin{array}{c|c}
 2& 1  \\
 1& 2  \\
\hline
 2&  
\end{array}
\end{smallmatrix}\right]
-
\left[\begin{smallmatrix}
\begin{array}{c|c}
 1&   \\
 1&   \\
\hline
  & 1  
\end{array}
\end{smallmatrix}\right]
=
\left[\begin{smallmatrix}
\begin{array}{c|c}
 1& 1  \\
  & 2  \\
\hline
 2&-1 
\end{array}
\end{smallmatrix}
\right],}\\
{M_{2,3}^t - I_{2,3}^t =
\left[\begin{smallmatrix}
\begin{array}{cc|c}
2 &  &1    \\
1 &1 &1    \\
1 &  &2    \\
\hline
1 &1 &     \\
  &2 &     
\end{array}
\end{smallmatrix}\right]
-
\left[\begin{smallmatrix}
\begin{array}{cc|c}
1 &  &     \\
1 &  &     \\
  &1 &     \\
\hline
  &  & 1   \\
  &  & 1   
\end{array}\end{smallmatrix}\right]
= 
\left[\begin{smallmatrix}
\begin{array}{cc|c}
 1 &  &1    \\
   &1 &1    \\
 1 &-1&2    \\
\hline
 1 & 1&-1   \\
   &2 &-1  
\end{array}
\end{smallmatrix}\right],}\\
{M_{3,4}^t - I_{3,4}^t =
\left[\begin{smallmatrix}
\begin{array}{ccc|cc}
2 &  &    &1   &    \\ 
1 &1 &    &1   &    \\
1 &  &1   &1   &    \\
  &1 &    &1   &1   \\
  &1 &    &    &2   \\
\hline
1 &  &1   &    &    \\
  &1 &1   &    &    \\
  &  &2   &    &    
\end{array}\end{smallmatrix}\right]
-
\left[\begin{smallmatrix}
\begin{array}{ccc|cc}
1 &  &    &  &   \\
1 &  &    &  &  \\
  &1 &    &  &  \\
  &  & 1  &  &   \\
  &  & 1  &  &   \\
\hline 
  &  &    &1 &   \\
  &  &    &1 &   \\
  &  &    &  & 1 
 \end{array}
 \end{smallmatrix}\right] 
= 
\left[\begin{smallmatrix}
\begin{array}{ccc|cc}
1 &  &    &1   &    \\ 
  &1 &    &1   &    \\
1 &-1&1   &1   &    \\
  &1 &-1  &1   &1   \\
  &1 &-1  &    &2   \\
\hline
1 &  &1   &-1  &    \\
  &1 &1   &-1  &    \\
  &  &2   &    &-1  
\end{array}\end{smallmatrix}\right],}
\end{gather*}}
and
{\allowdisplaybreaks
\begin{align*}
M_{4,5}^t - I_{4,5}^t 
&=
\left[\begin{smallmatrix}
\begin{array}{ccccc|ccc}
2 &  &    &    & &1  &  &  \\ 
1 &1 &    &    & &1  &  &  \\
1 &  &1   &    & &1  &  &  \\
  &1 &    &1   & &   &1 &  \\
  &1 &    &    &1&   &1 &  \\
  &  &1   &    & & 1 &  & 1\\
  &  &1   &    & &   &1 & 1\\
  &  &1   &    & &   &  & 2\\
\hline
1 &  &    & 1  & &   &  &  \\ 
  &1 &    & 1  & &   &  &  \\
  &  & 1  & 1  & &   &  &  \\
  &  &    & 1  &1&   &  &  \\
  &  &    &    &2&   &  &  
 \end{array}\end{smallmatrix}\right]
-
\left[\begin{smallmatrix}
\begin{array}{ccccc|ccc}
1 &  &    &  & &    &  &  \\
1 &  &    &  & &    &  & \\
  &1 &    &  & &    &  & \\
  &  & 1  &  & &    &  &  \\
  &  & 1  &  & &    &  &  \\
  &  &    &1 & &    &  &  \\
  &  &    &1 & &    &  &  \\
  &  &    &  &1&    &  &  \\
\hline
  &  &    &  & &1   &  &  \\
  &  &    &  & &1   &  & \\
  &  &    &  & &    &1 & \\
  &  &    &  & &    &  &1 \\
  &  &    &  & &    &  &1 
\end{array}\end{smallmatrix}\right] \\
& = 
\left[\begin{smallmatrix}
\begin{array}{ccccc|ccc}
1 &  &    &    &  &1  &  &  \\ 
  &1 &    &    &  &1  &  &  \\
1 &-1&1   &    &  &1  &  &  \\
  &1 &-1  &1   &  &   &1 &  \\
  &1 &-1  &    &1 &   &1 &  \\
  &  &1   &-1  &  & 1 &  & 1\\
  &  &1   & -1 &  &   &1 & 1\\
  &  &1   &    &-1&   &  & 2\\
\hline
1 &  &    & 1  &  &-1 &  &  \\ 
  &1 &    & 1  &  &-1 &  &  \\
  &  & 1  & 1  &  &   &-1&  \\
  &  &    & 1  &1 &   &  &-1\\
  &  &    &    &2 &   &  &-1
\end{array}\end{smallmatrix}\right].
\end{align*}}

It is easy to see that the kernels of the matrices 
$M_{l,l+1}^t - I_{l,l+1}^t$
are $\{ 0\}$ for all $l\in \Bbb N$.
Hence $K_1(\OLF) = \{0\}$ is obvious.
The computation of the $K_0$-group 
$K_0(\OLF)$ is the main body of this section.
We denote by ${\Bbb A}_{l+1,l}$
the $m(l+1) \times m(l)$  matrix 
$
 M_{l,l+1}^t - I_{l,l+1}^t.
 $
We will compute the cokernels of the matrices ${\Bbb A}_{l+1,l}$.
We set  subblock matrices 
${\Bbb A}_{l+1,l}^{UL},
{\Bbb A}_{l+1,l}^{UR},
{\Bbb A}_{l+1,l}^{LL}
$
and
${\Bbb A}_{l+1,l}^{LR}$
of ${\Bbb A}_{l+1,l}$ by setting
\begin{align*}
{\Bbb A}_{l+1,l}^{UL}(i,j) & = {\Bbb A}_{l+1,l}(i,j) \quad \text{ for } 
1 \le i \le m(l), 1\le j \le m(l-1),\\
{\Bbb A}_{l+1,l}^{UR}(i,j) & = {\Bbb A}_{l+1,l}(i,m(l-1)+j) \quad \text{ for } 
1 \le i \le m(l), 1\le j \le m(l-2),\\
{\Bbb A}_{l+1,l}^{LL}(i,j) & = {\Bbb A}_{l+1,l}(m(l) +i,j) \quad \text{ for } 
1 \le i \le m(l-1), 1\le j \le m(l-1),\\
{\Bbb A}_{l+1,l}^{LR}(i,j) & = {\Bbb A}_{l+1,l}(m(l) + i,m(l-1) + j) \quad \text{ for } 
1 \le i \le m(l-1), 1\le j \le m(l-2).
\end{align*}
They are 
an $m(l) \times m(l-1)$ matrix,  
an $m(l) \times m(l-2)$ matrix,
an $m(l-1) \times m(l-1)$ matrix and 
an $m(l-1) \times m(l-2) $ matrix 
respectively
such that
$$
{\Bbb A}_{l+1,l}
=
\left[
\begin{array}{c|c}
{\Bbb A}_{l+1,l}^{UL} & {\Bbb A}_{l+1,l}^{UR}\\
\hline
{\Bbb A}_{l+1,l}^{LL} & {\Bbb A}_{l+1,l}^{LR}\\
\end{array}
\right].
$$
Let $I_l$ be the $m(l)\times m(l)$ identity matrix. 
Recall that $0_{k,l}$ denotes
the $m(k)\times m(l)$ matrix all of which entries are $0's$.
By Lemma 3.1 and Lemma 3.2, one sees the general form of 
${\Bbb A}_{l+1,l}$ as in the following way:
\begin{lem}
 For $l=3,4,\dots,$, we have
\begin{align*}
 {\Bbb A}_{l+2,l+1}^{UL}
& = 
{\left[
\begin{array}{c|c|c}
\multicolumn{2}{c|}{} &  0_{l-1,l-2} \\
\multicolumn{2}{c|}{{\Bbb A}_{l+1,l}^{UL}}& \hrulefill \\
\multicolumn{2}{c|}{} &  I_{l-2} \\
\hline
 0_{l-1,l-2} &  S_{l-2}^t(1) &  {\Bbb A}_{l+1,l}^{LR}
\end{array}
\right],}\\
{\Bbb A}_{l+2,l+1}^{UR}
& =
{\left[
\begin{array}{c|c}
 S_{l-1}^t(1) &  0_{l,l-3}  \\
\hline
\multicolumn{2}{c}{ {\Bbb A}_{l+1,l}^{LL}}
\end{array}
\right],}\\
{\Bbb A}_{l+2,l+1}^{LL}
& = 
{\left[
\begin{array}{c|c}
I_{l-1}  &   \\
\hrulefill   & {\Bbb A}_{l+1,l}^{UR}\\
0_{l-2,l-1}  &  
\end{array}
\right],} \\
{\Bbb A}_{l+2,l+1}^{LR}
& =
{\left[
\begin{array}{c|c}
 {\Bbb A}_{l+1,l}^{LR}  & 0_{l-1,l-3}  \\
\hline
 0_{l-2,l-2}   & {\Bbb A}_{l,l-1}^{LR}  
\end{array}
\right].} 
\end{align*}
Hence the sequence ${\Bbb A}_{l+1,l}, l\in \Bbb N$  of the matrices
are
inductively determined.
\end{lem}

We set 
the 
$m(l) \times m(l)$ square matrix
$B_{l}$ by setting
$$
    B_{l} 
    = \left[
    \begin{array}{c|c}
{\Bbb A}_{l+1,l}^{UL} & {\Bbb A}_{l+1,l}^{UR}
\end{array}
\right]
$$
the upper  half of the matrix ${\Bbb A}_{l+1,l}$.
We next provide a sequence 
 $C_{l+1,l}, l\in \Bbb N $ 
 of $m(l+1) \times m(l)$ matrix such as:

$$
C_{2,1}=
\left[
\begin{array}{c|c}
 2&    \\
 2&   \\
\hline
 2&  
\end{array}
\right],
 \qquad
C_{3,2}              
\left[
\begin{array}{c|c|c}
1 &1 &     \\
1 &1 &     \\
1 &1 &     \\
\hline
1 &1 &     \\
1 &1 &     \\
\end{array}
\right],
\qquad
C_{4,3} =
\left[
\begin{array}{cc|c|cc}
1 &  &1  &  &    \\ 
1 &  &1  &  &    \\
1 &  &1  &  &    \\
1 &  &1  &  &    \\
1 &  &1  &  &    \\
\hline
  &1 &1  &  &    \\
  &1 &1  &  &    \\
  &1 &1  &  &    \\
\end{array}
\right],
$$
$$
C_{5,4} =
\left[
\begin{array}{ccc|cc|ccc}
1 &  &   &1 &  &   &  &  \\ 
1 &  &   &1 &  &   &  &  \\
1 &  &   &1 &  &   &  &  \\
1 &  &   &1 &  &   &  &  \\
1 &  &   &1 &  &   &  &  \\
  &1 &   &1 &  &   &  &  \\
  &1 &   &1 &  &   &  &  \\
  &1 &   &1 &  &   &  &  \\
\hline
  &  & 1 &  &1&   &  &  \\ 
  &  & 1 &  &1&   &  &  \\
  &  & 1 &  &1&   &  &  \\
  &  & 1 &  &1&   &  &  \\
  &  & 1 &  &1&   &  &  
 \end{array}
\right],
\qquad
C_{6,5} =
\left[
\begin{array}{ccccc|ccc|ccccc}
1 &  &  & &   &1 & &  & & &   &  &  \\ 
1 &  &  & &   &1 & &  & & &   &  &  \\
1 &  &  & &   &1 & &  & & &   &  &  \\
1 &  &  & &   &1 & &  & & &   &  &  \\
1 &  &  & &   &1 & &  & & &   &  &  \\
  & 1&  & &   &1 & &  & & &   &  &  \\
  & 1&  & &   &1 & &  & & &   &  &  \\
  & 1&  & &   &1 & &  & & &   &  &  \\
  &  & 1& &   &  &1&  & & &   &  &  \\
  &  & 1& &   &  &1&  & & &   &  &  \\
  &  & 1& &   &  &1&  & & &   &  &  \\
  &  & 1& &   &  &1&  & & &   &  &  \\
  &  & 1& &   &  &1&  & & &   &  &  \\
\hline
  &  &  &1&   &  & &1 & & &   &  &  \\ 
  &  &  &1&   &  & &1 & & &   &  &  \\
  &  &  &1&   &  & &1 & & &   &  &  \\
  &  &  &1&   &  & &1 & & &   &  &  \\
  &  &  &1&   &  & &1 & & &   &  &  \\
  &  &  & & 1 &  & &1 & & &   &  &  \\
  &  &  & & 1 &  & &1 & & &   &  &  \\
  &  &  & & 1 &  & &1 & & &   &  &  
 \end{array}
 \right].
$$
To define the matrices $C_{l+1,l}$ for $l \ge 6$, 
divide $C_{l+1,l}$ into $6$ subblock matrices
 $C_{l+1,l}^{UL},C_{l+1,l}^{UM},C_{l+1,l}^{UR},C_{l+1,l}^{LL},C_{l+1,l}^{LM},C_{l+1,l}^{LR}$ as in the following way: 
\begin{align*}
C_{l+1,l}^{UL}(i,j) = & C_{l+1,l}(i,j) \text{ for } 1\le i \le m(l),
1\le j\le m(l-2), \\  
C_{l+1,l}^{UM} (i,j) = & C_{l+1,l}(i,j+m(l-2)) \text{ for } 1\le i \le m(l),
1\le j\le m(l-3), \\
C_{l+1,l}^{UR} (i,j) = & C_{l+1,l}(i,j+m(l-2)+m(l-3)) \text{ for } 1\le i \le m(l),
1\le j\le m(l-2), \\
C_{l+1,l}^{LL}(i,j) = & C_{l+1,l}(i+m(l),j) \text{ for } 1\le i \le m(l-1),
1\le j\le m(l-2), \\  
C_{l+1,l}^{LM} (i,j) = & C_{l+1,l}(i+m(l),j+m(l-2)) \text{ for } 1\le i \le m(l-1),
1\le j\le m(l-3), \\
C_{l+1,l}^{LR} (i,j) = & C_{l+1,l}(i+m(l),j+m(l-2)+m(l-3)) \text{ for } 1\le i \le m(l-1),
1\le j\le m(l-2). 
\end{align*}
They are 
an $ m(l) \times m(l-2)$ matrix,  
 an $ m(l) \times m(l-3)$ matrix,
 an $ m(l) \times m(l-2)$ matrix,
 an $ m(l-1) \times m(l-2)$ matrix,
 an $m(l-1) \times m(l-3)$ matrix and 
 an 
$m(l-1) \times m(l-2) $ matrix 
respectively
such that

$$
C_{l+1,l}
=
\left[
\begin{array}{c|c|c}
C_{l+1,l}^{UL} & C_{l+1,l}^{UM} & C_{l+1,l}^{UR}\\
\hline
C_{l+1,l}^{LL} & C_{l+1,l}^{LM} & C_{l+1,l}^{LR}\\
\end{array}
\right].
$$
These block matrices are defined inductively as in the following way:

\begin{gather*}
{C_{l+1,l}^{UL}
= 
\left[
\begin{array}{c|c}
C_{l,l-1}^{UL}  &   \\
\hrulefill         &  0_{l,l-4}\\
C_{l,l-1}^{LL}  &   
\end{array}
\right], \quad
C_{l+1,l}^{UM}
 = 
\left[
\begin{array}{c|c}
C_{l,l-1}^{UM}  &   \\
\hrulefill      & 0_{l,l-5}\\
C_{l,l-1}^{LM}  &  
\end{array}
\right], \quad
C_{l+1,l}^{UL}
 =
 \left[
0_{l,l-2}
\right], }\\
{C_{l+1,l}^{LL}
 = 
\left[
\begin{array}{c|c}
                & C_{l-1,l-2}^{UL} \\
0_{l-1,l-3}     & \hrulefill \\
                & C_{l-1,l-2}^{LL} 
\end{array}
\right], \quad
C_{l+1,l}^{LM}
 = 
\left[
\begin{array}{c|c}
                & C_{l-1,l-2}^{UM} \\
0_{l-1,l-4}     & \hrulefill \\
                & C_{l-1,l-2}^{LM} 
\end{array}
\right], \quad
C_{l+1,l}^{LR}
 =
\left[
0_{l,l-2}
\right].}
 \end{gather*}
Let $L_{l+1,l}$ be the $m(l+1)\times m(l)$ matrix defined by the 
block matrix:
$$
L_{l+1,l}
=
\left[
\begin{array}{c|c}
L_{l+1,l}^{UL} & L_{l+1,l}^{UR}\\
\hline
L_{l+1,l}^{LL} & L_{l+1,l}^{LR}
\end{array}
\right]
$$
where
\begin{align*}
L_{l+1,l}^{UL}
&  = {\Bbb A}_{l+1,l}^{UL}: \quad m(l) \times m(l-1)
 \text{ matrix, } \\
L_{l+1,l}^{UR}
&  =
\left[\begin{smallmatrix}
& 0_{l-1,l-2}&  \\
\hline
& B_{l-2}& 
\end{smallmatrix}\right] :\quad m(l) \times m(l-2)
 \text{ matrix, }\\
L_{l+1,l}^{LL}
&  ={\Bbb A}_{l+1,l}^{LL}: \quad m(l-1) \times m(l-1)\text{ matrix,} \\
L_{l+1,l}^{LR} 
& = -I^{t}_{l-2,l-1} -C_{l-1,l-2}: \quad m(l-1) \times m(l-2)
 \text{ matrix}.
\end{align*}
We write down the above matrices for $l=1,2,3,4.$
$$
L_{2,1}
=
\left[
\begin{array}{c|c}
 1&   \\
  &2  \\
\hline
 2&-3 \\
\end{array}
\right],
\qquad
L_{3,2}= 
\left[
\begin{array}{cc|c}
 1 &  &     \\
   &1 &     \\
 1 &-1&2    \\
\hline
 1 & 1&-3   \\
   &2 &-3   \\
\end{array}
\right],
\qquad
L_{4,3}= 
\left[
\begin{array}{ccc|cc}
1 &  &    &    &    \\ 
  &1 &    &    &    \\
1 &-1&1   &    &    \\
  &1 &-1  &1   &1   \\
  &1 &-1  &    &2   \\
\hline
1 &  &1   &-3  &    \\
  &1 &1   &-3  &    \\
  &  &2   &-2  &-1  
\end{array}
\right],
$$
and
$$
L_{5,4}= 
\left[
\begin{array}{ccccc|ccc}
1 &  &    &    &  &   &  &  \\ 
  &1 &    &    &  &   &  &  \\
1 &-1&1   &    &  &   &  &  \\
  &1 &-1  &1   &  &   &  &  \\
  &1 &-1  &    &1 &   &  &  \\
  &  &1   &-1  &  & 1 &  & 1\\
  &  &1   & -1 &  &   &1 & 1\\
  &  &1   &    &-1& 1 &-1& 2\\
\hline
1 &  &    & 1  &  &-2 &-1&  \\ 
  &1 &    & 1  &  &-2 &-1&  \\
  &  & 1  & 1  &  &-1 &-2&  \\
  &  &    & 1  &1 &-1 &-1&-1\\
  &  &    &    &2 &-1 &-1&-1
\end{array}
\right].
$$
We define the elementary column operations on integer matrices to be:
\begin{enumerate}
\item Multiply a column by $-1$,
\item Add an integer multiple of one column to another column.
\end{enumerate}

The elementary row operations are similarly defined.
We know that the matrices $L_{l+1,l}$ is obtained from ${\Bbb A}_{l+1,l}$
by elementary column operations, that operation is denoted by $\Gamma_l$.
The operation
$\Gamma_l$ is an $m(l) \times m(l)$ matrix corresponding to the column operation such that
$$
L_{l+1,l} = {\Bbb A}_{l+1,l} \Gamma_l.
$$
Since
$$
L_{l+1,l}
= 
\left[
\begin{array}{c|c}
 & 0_{l-1,l-2} \\
{\Bbb A}_{l+1,l}^{UL}  & \hrulefill \\
 & B_{l-2} \\
\hline
 {\Bbb A}_{l+1,l}^{LL} &  -I^t_{l-2,l-1} - C_{l-1,l-2}
\end{array}
\right],
$$
we may apply the elementary column operation 
$I_{l-1} \oplus \Gamma_{l-2}$
to
$L_{l+1,l}$
so that 
 the matrix
$B_{l-2}$ in $L_{l+1,l}$ 
goes to
$$
\left[
\begin{array}{c|c}
&  0_{l-3,l-1}                       \\
{\Bbb A}_{l-1,l-2}^{UL} & \hrulefill \\ 
  &  B_{l-4}\\
\end{array}
\right].
$$
The new matrix 
$L_{l+1,l}(I_{l-1} \oplus \Gamma_{l-2})$
 is
$$
L_{l+1,l}(I_{l-1} \oplus \Gamma_{l-2})
= 
\left[
\begin{array}{c|c|c}
&  \multicolumn{2}{c}{ 0_{l-1,l-2} }\\
{\Bbb A}_{l+1,l}^{UL} & \multicolumn{2}{c}{\hrulefill} \\
              &  & 0_{l-3,l-4} \\
              & {\Bbb A}_{l-1,l-2}^{UL} & \hrulefill \\
              &  &     B_{l-4} \\
\hline
 {\Bbb A}_{l+1,l}^{LL} & 
 \multicolumn{2}{c}{ (-I^t_{l-2,l-1} - C_{l-1,l-2})\Gamma_{l-2}}
\end{array}
\right].
$$
As 
$$
B_{l-2n} \Gamma_{l-2n} =
\left[
\begin{array}{c|c}
&  0_{l-2n-1,l-2n-2}\\
{\Bbb A}_{l-2n+1,l-2n}^{UL}& \hrulefill \\ 
& B_{l-2n-2}\\
\end{array}
\right]
$$
for $n=1,2,\dots$ with $2n <l $,
by continuing these procedures $k$-times for $l  = 2k, 2k+1$
we finally get 
$$
B_{2}\Gamma_2 = 
 \begin{bmatrix}
 1 & 0 &0\\
0 & 1 &0\\
1 & -1& 2
\end{bmatrix}
\quad \text{ for } l=2k
\quad \text{ and }
\quad
B_{1}\Gamma_1 = 
 \begin{bmatrix}
1 &0 \\ 
0& 2
\end{bmatrix}
\quad \text{ for } l=2k+1.
$$
For $l=2k, 2k+1$,
let 
${\Bbb M}_{l+1,l}$ be the $m(l+1) \times m(l)$ 
matrix obtained from $L_{l+1,l}$
after the $k$ times procedures above. 
Then we have
$$
{\Bbb M}_{l+1,l}(i,j) =
\begin{cases}
0 & \text{ if } i < j, \, 1 \le i,j \le m(l) \\
1 & \text{ if } i = j, \, 1 \le i < m(l) \\ 
2 & \text{ if } i = j= m(l).
\end{cases}
$$ 
Let
$
v_l = 
[v_l(i)]_{i=1}^{m(l-1)}
$
be the column vector of length $m(l-1)$ defined by
$$
v_l(i) = {\Bbb M}_{l+1,l}(m(l) + i,m(l)), \qquad i=1,2,\dots, m(l-1)
$$
so that the matrix ${\Bbb M}_{l+1,l}$ is of the form
$$
{\Bbb M}_{l+1,l} = 
\begin{bmatrix}
1 & & & &   \\
  &1& & &   \\
  & &\ddots   & & \\
  &* & &1&   \\
  & & & &2  \\
 \hline
  & & & & v_l(1)   \\
  & & & & v_l(2)  \\
  &* & & & \vdots  \\
  & & & & v_l(m(l-1))  
\end{bmatrix}.
$$
For 
$l=1,2,3,4,5,6,$ we see
$$
v_1
= 
\left[-3
\right],
\quad
v_2
= 
\left[\begin{smallmatrix}
-3\\
-3\\
\end{smallmatrix}\right],
\quad
v_3
= 
\left[\begin{smallmatrix}
3\\
3\\
1 \\
\end{smallmatrix}
\right],
\quad
v_4
= 
\left[\begin{smallmatrix}
3\\
3\\
3\\
1 \\
1 \\
\end{smallmatrix}
\right],
\quad
v_5
= 
\left[\begin{smallmatrix}
-3\\
-3\\
-3\\
-3\\
-3\\
-1 \\
-1 \\
-3\\
\end{smallmatrix}
\right],
\quad
v_6
= 
\left[\begin{smallmatrix}
-3\\
-3\\
-3\\
-3\\
-3\\
-3\\
-3\\
-3\\
-1 \\
-1 \\
-1 \\
-3\\
-3\\
\end{smallmatrix}
\right].
$$
By induction, one has:
\begin{lem}\hspace{6cm}
\begin{enumerate}
\renewcommand{\labelenumi}{(\roman{enumi})}
\item 
$
v_l(i) = 
\begin{cases}
-3 & \text{ if } l= 4k+1,4k+2, k\in \Zp, \text{ and } 1 \le i \le m(l-2),\\ 
3 & \text{ if } l= 4k+3,4k+4, k\in \Zp, \text{ and } 1 \le i \le m(l-2),\\
\end{cases}
$
\item
$v_l(m(l-2) + i) = \widehat{v_{l-2}(i)}$ for $ i=1,2,\dots, m(l-3)$
\end{enumerate}
where
for $u=\pm 3, \pm 1,$ the integer $\hat{u}$ is defined by
$$
\hat{u}
=
\begin{cases}
u-4 & \text{ if } u = 3,1, \\
u+4 & \text{ if } u =- 3,-1. \\
\end{cases}
$$ 
\end{lem}
\begin{align*}
{\Bbb N}_{l+1,l}(i,j)
& =
\begin{cases} 
1 & \text{ if } i=j, \, 1 \le i < m(l), \\
2 & \text{ if } i=j= m(l), \\
v_l(i-m(l)) & \text{ if } i > m(l),\,  j= m(l), \\
0 & \text{ otherwise,}
\end{cases} \\
{\Bbb H}_{l+1,l}(i,j)
& =
\begin{cases} 
1 & \text{ if } i=j, \, 1 \le i < m(l), \\
2 & \text{ if } i=j= m(l), \\
-1 & \text{ if } i > m(l),\,  j= m(l), \\
0 & \text{ otherwise.}
\end{cases}
\end{align*}
We set 
$m(l+1) \times m(l)$ matrices 
${\Bbb N}_{l+1,l}$ and 
${\Bbb H}_{l+1,l}$ 
by setting
For $l=1,2,3,4,$ one sees 
$$
{\Bbb N}_{2,1}
=
\left[\begin{smallmatrix}
 1&   \\
  &2  \\
\hline
  &-3 \\
\end{smallmatrix}\right],
\quad
{\Bbb N}_{3,2}= 
\left[\begin{smallmatrix}
 1 &   &     \\
   &1  &     \\
   &   &2    \\
\hline
   &  &-3   \\
   &  &-3   \\
\end{smallmatrix}\right],
\quad
{\Bbb N}_{4,3}= 
\left[\begin{smallmatrix}
1 &  &    &    &    \\ 
  &1 &    &    &    \\
  &  &1   &    &    \\
  &  &    &1   &    \\
  &  &    &    &2   \\
\hline
  &  &    &    & 3     \\
  &  &    &    & 3     \\
  &  &    &    & 1  \\
\end{smallmatrix}
\right],
\quad
{\Bbb N}_{5,4}= 
\left[\begin{smallmatrix}
1 &  &    &    &  &   &  &  \\ 
  &1 &    &    &  &   &  &  \\
  &  &1   &    &  &   &  &  \\
  &  &    &1   &  &   &  &  \\
  &  &    &    &1 &   &  &  \\
  &  &    &    &  & 1 &  &  \\
  &  &    &    &  &   &1 &  \\
  &  &    &    &  &   &  & 2\\
\hline
  &  &    &    &  &   &  &3 \\ 
  &  &    &    &  &   &  &3 \\
  &  &    &    &  &   &  &3 \\
  &  &    &    &  &   &  &1\\
  &  &    &    &  &   &  &1\\
\end{smallmatrix}
\right]
$$
and
$$
{\Bbb H}_{2,1}
=
\left[\begin{smallmatrix}
 1&   \\
  &2  \\
\hline
  &-1 \\
\end{smallmatrix}\right],
\quad
{\Bbb H}_{3,2}= 
\left[\begin{smallmatrix}
 1 &   &     \\
   &1  &     \\
   &   &2    \\
\hline
   &  &-1   \\
   &  &-1   \\
\end{smallmatrix}\right],
\quad
{\Bbb H}_{4,3}= 
\left[\begin{smallmatrix}
1 &  &    &    &    \\ 
  &1 &    &    &    \\
  &  &1   &    &    \\
  &  &    &1   &    \\
  &  &    &    &2   \\
\hline
  &  &    &    &-1     \\
  &  &    &    &-1     \\
  &  &    &    &-1  \\
\end{smallmatrix}
\right],
\quad
{\Bbb H}_{5,4}= 
\left[\begin{smallmatrix}
1 &  &    &    &  &   &  &  \\ 
  &1 &    &    &  &   &  &  \\
  &  &1   &    &  &   &  &  \\
  &  &    &1   &  &   &  &  \\
  &  &    &    &1 &   &  &  \\
  &  &    &    &  & 1 &  &  \\
  &  &    &    &  &   &1 &  \\
  &  &    &    &  &   &  & 2\\
\hline
  &  &    &    &  &   &  &-1\\ 
  &  &    &    &  &   &  &-1\\
  &  &    &    &  &   &  &-1\\
  &  &    &    &  &   &  &-1\\
  &  &    &    &  &   &  &-1\\
\end{smallmatrix}
\right].
$$
By elementary row operations compatible to $I_{l,l+1}^t$,
one gets the matrix ${\Bbb N}_{l+1,l}$ 
from the matrix ${\Bbb M}_{l+1,l}$.
In the matrix ${\Bbb N}_{l+1,l}$,
for $i=1,2,\dots,m(l)$,
if $v_l(i) = -3,$ 
then add the $m(l)$-th row to the $i+m(l)$-th row at the $i+m(l)$-th row, 
if $v_l(i) = 3,$ 
then subtract the twice of $m(l)$-th row from the $i+m(l)$-th row at the $i+m(l)$-th row,
if $v_l(i) = -3,$ 
then subtract the $m(l)$-th row from the $i+m(l)$-th row 
at the $i+m(l)$-th row,
then one gets the matrix
${\Bbb H}_{l+1,l}$.
These row operations are compatible to the map
$I_{l,l+1}^t$
and the relations
$$
I^t_{l,l+1}{\Bbb N}_{l,l-1}= {\Bbb N}_{l+1,l}I^t_{l-1,l},
\qquad
I^t_{l,l+1}{\Bbb H}_{l,l-1}= {\Bbb H}_{l+1,l}I^t_{l-1,l} 
$$
for $l =2,3,\dots $ hold.
As 
$$
(M_{l,l+1}^t - I_{l,l+1}^t ){\Bbb Z}^{m(l)}
= {\Bbb A}_{l+1,l}{\Bbb Z}^{m(l)} 
= L_{l+1,l}{\Bbb Z}^{m(l)} 
= {\Bbb M}_{l+1,l}{\Bbb Z}^{m(l)}, \qquad l \in {\Bbb N}
$$
we see that  
$
{\Bbb Z}^{m(l+1)}/(M_{l,l+1}^t - I_{l,l+1}^t ){\Bbb Z}^{m(l)}
$
coincides with the group
$
{\Bbb Z}^{m(l+1)} /{\Bbb M}_{l+1,l}{\Bbb Z}^{m(l)}
$
 for all $l \in \Bbb N$.
We then have
\begin{prop}
There exist isomorphisms
\begin{align*}
\xi_l & : {\Bbb Z}^{m(l)}/{\Bbb M}_{l,l-1}{\Bbb Z}^{m(l-1)} \rightarrow
{\Bbb Z}^{m(l)}/{\Bbb N}_{l,l-1}{\Bbb Z}^{m(l-1)},\\
\eta_l & : {\Bbb Z}^{m(l)}/{\Bbb N}_{l,l-1}{\Bbb Z}^{m(l-1)} \rightarrow
{\Bbb Z}^{m(l)}/{\Bbb H}_{l,l-1}{\Bbb Z}^{m(l-1)}
\end{align*}
of abelian groups such that
the following diagrams are commutative:
$$
\begin{CD}
{\Bbb Z}^{m(l)}/(M_{l-1,l}^t - I_{l-1,l}^t ){\Bbb Z}^{m(l-1)}
@>
\bar{I}_{l,l+1}^t 
>> 
{\Bbb Z}^{m(l+1)}/(M_{l,l+1}^t - I_{l,l+1}^t ){\Bbb Z}^{m(l)}
 \\
@| @| \\
{\Bbb Z}^{m(l)}/{\Bbb M}_{l,l-1}{\Bbb Z}^{m(l-1)}
@. 
{\Bbb Z}^{m(l+1)}/{\Bbb M}_{l+1,l}{\Bbb Z}^{m(l)}
 \\
@V
\xi_l  
VV @V
\xi_{l+1} 
VV \\
{\Bbb Z}^{m(l)}/ {\Bbb N}_{l,l-1}{\Bbb Z}^{m(l)}
@. 
{\Bbb Z}^{m(l+1)}/ {\Bbb N}_{l+1,l}{\Bbb Z}^{m(l+1)}\\
@V
\eta_l  
VV @V
\eta_{l+1} 
VV \\
{\Bbb Z}^{m(l)}/ {\Bbb H}_{l,l-1}{\Bbb Z}^{m(l)}
@>\widehat{I}_{l,l+1}^t >> 
{\Bbb Z}^{m(l+1)}/ {\Bbb H}_{l+1,l}{\Bbb Z}^{m(l+1)}
\end{CD}
$$
where
$\widehat{I}_{l,l+1}^t:
{\Bbb Z}^{m(l)}/ {\Bbb H}_{l,l-1}{\Bbb Z}^{m(l-1)}
\rightarrow 
{\Bbb Z}^{m(l+1)}/ {\Bbb H}_{l+1,l}{\Bbb Z}^{m(l)}
$
is the homomorphism induced by the matrix
$I_{l,l+1}^t.$
Hence
we have an isomorphism
$$
K_0(\OLF)
 \cong
 \underset{l}{\varinjlim} 
\{ \widehat{I}_{l,l+1}^t:
{\Bbb Z}^{m(l)}/{\Bbb H}_{l,l-1}{\Bbb Z}^{m(l-1)}
\rightarrow
{\Bbb Z}^{m(l+1)}/{\Bbb H}_{l+1,l}{\Bbb Z}^{m(l)}
 \}.
$$
\end{prop}
 We fix $l \ge 3$. 
 Define the $(m(l-1)+1) \times 1$ matrix $R_{l-1}$ and 
 the $ (m(l-1) +1) \times (m(l-2) +1)$ matrix 
 $I_{l-1,l-2}^R$ by setting:
$$
  R_{l-1} =
\left[\begin{smallmatrix}
2\\
-1\\ 
\vdots\\
-1\\
\end{smallmatrix}\right],
\qquad
I_{l-1,l-2}^R
=
\left[
\begin{array}{c|ccc}
1       & 0 &\dots &0 \\
\hline
0       & &      & \\
\vdots  & &   I_{l-2,l-1} & \\
0       & &      & 
\end{array}\right].
$$
Then the following diagram is commutative:
$$
\begin{CD}
{\Bbb Z}^{m(l)}/{\Bbb H}_{l,l-1}{\Bbb Z}^{m(l-1)}
@>
\widehat{I}_{l,l+1}^t 
>> 
{\Bbb Z}^{m(l+1)}/{\Bbb H}_{l+1,l}{\Bbb Z}^{m(l)}
 \\
@| @| \\
{\Bbb Z}^{m(l-2)+1}/{R}_{l-2}{\Bbb Z}
@>
\bar{I}_{l-1,l-2}^R 
>> 
{\Bbb Z}^{m(l-1)+1}/ R_{l-1}{\Bbb Z}
 \end{CD}
$$
where $\bar{I}_{l-1,l}^R$ is the homomorphism induced by the matrix
$I_{l-1,l}^R.$
Let
$\varphi_{l-2}:{\Bbb Z}^{m(l-2) +1} \rightarrow 
{\Bbb Z}^{m(l-2) +1} $
be an isomorphism
defined by the operations on the row vectors of ${\Bbb Z}^{m(l-2) +1} $ 
to 
add the 2-times multiplication of the second row to the first row,
and subtract the second row from the $k$-th rows for
$k=3,4,\dots, m(l-2) +1$.
It is given by the matrix:
$$
Q_{l-2}
=
\left[\begin{smallmatrix}
1 & 2       &     &  &       &      \\
  & 1       &     &  &       &      \\
  &-1       &1    &  &       &      \\
  &-1       &     &1 &       &      \\
  & \vdots  &     &  &\ddots &      \\
  &-1       &     &  &       & 1    \\
\end{smallmatrix}\right].
$$ 
 Since
 $
 Q_{l-2} R_{l-2}
 =
 \left[\begin{smallmatrix}
 0\\
 -1\\
0\\
\vdots\\
0
\end{smallmatrix}\right],
$
$\varphi_{l-2}$ yields an isomorphism
$$
\varphi_{l-2}: {\Bbb Z}^{m(l-2)+1}/R_{l-2}{\Bbb Z}
\rightarrow
{\Bbb Z} \oplus 0 \oplus  {\Bbb Z}^{m(l-2)-1} = {\Bbb Z}^{m(l-2)}.
$$
Let $J_{l-1,l-2}: {\Bbb Z}^{m(l-2)-1} \rightarrow {\Bbb Z}^{m(l-1)-1}$ 
be a homomorphism defined by the $(m(l-1)-1) \times (m(l-2)-1)$ matrix
$$
J_{l-1,l-2}(i,j) =
\begin{cases}
0 & \text{ if } i=1, \\
I_{l-2,l-1}(i+1, j+1) & \text{ if } i = 2,\dots,m(l-2)-1
\end{cases}
$$
for $i=1,2,\dots,m(l-1)-1,\, j=1,2,\dots, m(l-2)-1$.
We set 
$\widetilde{I}_{l-1,l-2}: {\Bbb Z}^{m(l-2)} \rightarrow {\Bbb Z}^{m(l-1)}$ 
a homomorphism defined by the $m(l-1) \times m(l-2)$ matrix 
$$
\widetilde{I}_{l-1,l-2}(i,j) =
\begin{cases}
1 & \text{ if } i=j = 1, \\
0 & \text{ if } i=1, \,  j \ge 2, \\
0 & \text{ if } i=2, \\
I_{l-2,l-1}(i, j) & \text{ if } i = 3,4,\dots,m(l-2)-1
\end{cases}
$$
for $i=1,2,\dots,m(l-1),\, j=1,2,\dots, m(l-2)$.
That is,
$$
\widetilde{I}_{l-1,l-2}
=
\left[\begin{smallmatrix}
\begin{array}{c|ccc}
1       & 0 &\dots &0 \\
\hline
0       &   &      & \\
\vdots  &   & J_{l-1,l-2} & \\
0       &   &          & 
\end{array}\end{smallmatrix}\right].
$$
\begin{lem}
The diagram 
$$
\begin{CD}
{\Bbb Z}^{m(l-2)+1}/{R}_{l-2}{\Bbb Z}
@>
\bar{I}^ R_{l-1,l-2} 
>> 
{\Bbb Z}^{m(l-1)+1}/ R_{l-1}{\Bbb Z}
 \\
@V
\varphi_{l-2}  
VV @V
\varphi_{l-1} 
VV \\
 {\Bbb Z}^{m(l-2)}
@>\widetilde{I}_{l-1,l-2} >> 
 {\Bbb Z}^{m(l-1)}
\end{CD}
$$
is  commutative. 
Hence
we have an isomorphism
$$
K_0(\OLF)
\cong
{\Bbb Z}
\oplus
{ \underset{l}{\varinjlim} 
\{ 
J_{l-1,l-2}:
{\Bbb Z}^{m(l-2)-1}
\rightarrow
{\Bbb Z}^{m(l-1)-1} \}.}
$$
\end{lem}
\begin{pf}
Since 
the commutativity
$\varphi_{l-1}\circ \bar{I}^ R_{l-1,l-2} 
= \widetilde{I}_{l-1,l-2}\circ \varphi_{l-2}
$
is immediate, one has 
$$
K_0(\OLF)
\cong
 \underset{l}{\varinjlim} 
\{ \widetilde{I}_{l-1,l-2}:
{\Bbb Z}^{m(l-2)}
\rightarrow
{\Bbb Z}^{m(l-1)} \}.
$$
As 
$\widetilde{I}_{l-1,l-2} = 1 \oplus J_{l-1,l-2}$,
the assertion is clear.
\end{pf}
We will compute the group of the inductive limit
$ \underset{l}{\varinjlim} 
\{ J_{l+1,l}:
{\Bbb Z}^{m(l)-1}
\rightarrow
{\Bbb Z}^{m(l+1)-1}\}
$,
that we denote by $G$.
Let $I^c_{l+1,l}$ be the $(m(l)-2 )\times (m(l) -1) $ matrix 
defined by 
$$
I^c_{l+1,l}(i,j) = I_{l,l+1}(i+2,j+1)
\quad \text{ for} \quad i=1,\dots,m(l)-2,
\, \,
 j=1,\dots,m(l)-1.
$$
Hence
$
J_{l+1,l}
= 
\left[\begin{smallmatrix}
0 &   \cdots &   0 \\
\hline
  &               &   \\
  &      I_{l+1,l}^c   &   \\
  &            &   \\
\end{smallmatrix}\right].
$
It gives rise to a homomorphism :
$$
I^c_{l+1,l}: {\Bbb Z}^{m(l)-1} \rightarrow 
{\Bbb Z}\oplus I^c_{l+1,l}{\Bbb Z}^{m(l)-1} \subset {\Bbb Z}^{m(l)-1}.
$$
Put
$$
{\Bbb Z}(l) = {\Bbb Z} \oplus \cdots \oplus {\Bbb Z} = {\Bbb Z}^{m(l)-1}.
$$
For $k \in \Bbb N$, take $l \in \Zp$ such that 
$k \le m(l)$.
Define a sequence of positive integers
$$
g_k = \sum_{j=1}^k \sum_{i=2}^{m(l+1)} I_{l,l+1}^t(i,j), \qquad k=1,2,\dots
$$
that is independent of the choice of $l$, 
so that
$$
g_1 = 1,\quad
g_2 = 2,\quad
g_3 = 4,\quad
g_4 = 6,\quad
g_5 = 7,\quad
g_6 = 9,
\dots.
$$
Define for $l \ge k$,
$$
{\Bbb Z}(l;k)
= 
\overbrace{0\oplus \cdots \oplus 0}^{g_k } 
\oplus 
\overbrace{{\Bbb Z}\oplus \cdots \oplus{\Bbb Z}}^{ m(l)-1-g_k }
\subset
{\Bbb Z}^{m(l)-1}= {\Bbb Z}(l)
$$
so that 
we have
$$
I_{l+1,l}^c({\Bbb Z}(l;k)) \subset {\Bbb Z}(l+1;k+1).
$$
Set the group of the inductive limit
$$
G_k =
\underset{n}{\varinjlim} 
\{ I_{k+n+1,k+n}^c : {\Bbb Z}(k+n;n)\rightarrow {\Bbb Z}(k+n+1;n+1) \}.
$$
Since the following diagram is commutative:
$$
\begin{CD}
{\Bbb Z}(1)
@>
I_{2,1}^c 
>> 
{\Bbb Z}(2;1)
@>
I_{3,2}^c 
>> 
{\Bbb Z}(3;2)
@>
I_{4,3}^c 
>> 
{\Bbb Z}(4;3)
@>
I_{5,4}^c 
>> 
\cdots 
@>>>
G_1 \\
@.@VV{\iota}V @VV{\iota}V @VV{\iota}V @. \\
@.
{\Bbb Z}(2)
@>
I_{3,2}^c 
>> 
{\Bbb Z}(3;1)
@>
I_{4,3}^c 
>> 
{\Bbb Z}(4;2)
@>
I_{5,4}^c 
>> 
\cdots  
@>>> G_2 \\
@. @. @VV{\iota}V @VV{\iota}V @. \\
@. @.
{\Bbb Z}(4)
@>
I_{4,3}^c 
>> 
{\Bbb Z}(4;1)
@>
I_{5,4}^c 
>> 
\cdots  
@>>> G_3 \\
@. @. @. @VV{\iota}V @. \\
@. @. @.
{\Bbb Z}(4)
@>
I_{5,4}^c 
>> 
\cdots  
@>>> G_4 \\
@. @. @.
@. 
\vdots  
\end{CD}
$$
where
the vertical arrows $\iota$ mean the natural inclusion maps,
one sees the next lemma:
\begin{lem}\hspace{8cm}
\begin{enumerate}
\renewcommand{\labelenumi}{(\roman{enumi})}
\item  For each $k =1,2,\dots$, 
the group $G_k$ is isomorphic to the abelian group 
$C({\frak K}_k,{\Bbb Z})$ of all integer valued continuous functions on a Cantor discontinuum ${\frak K}_k$.
\item The sequence $G_k, k=1,2,\dots $ are increasing whose union generate
$G$.
\end{enumerate}
\end{lem}
Hence one has
\begin{lem}
The group $G$ is isomorphic to the countable direct sum of the group 
$C({\frak K},{\Bbb Z})$ of all integer valued continuous functions on a Cantor discontinuum ${\frak K}$.
\end{lem}
\begin{pf}
It is easy to see that
$G_k$ is isomorphic to the direct sum 
$C({\frak K}_{k,k-1},{\Bbb Z}) \oplus G_{k-1}$
of all integer valued continuous functions on a Cantor discontinuum 
${\frak K}_{k,k-1}$ and $G_{k-1}$ for
each $k$.
Hence we have
\begin{align*}
G_k
& \cong C({\frak K}_{k,k-1},{\Bbb Z}) \oplus G_{k-1}\\ 
& \cong C({\frak K}_{k,k-1},{\Bbb Z}) \oplus C({\frak K}_{k-1,k-2},{\Bbb Z}) 
\oplus \cdots \oplus C({\frak K}_{2,1},{\Bbb Z}) \oplus G_{1}.
\end{align*}
Since both $G_1$ and  $C({\frak K}_{i,i-1},{\Bbb Z})$
are isomorphic to the group $C({\frak K},{\Bbb Z})$
of all integer valued continuous functions on a Cantor discontinuum
$\frak K$, 
we have
$$
G \cong \underset{k}{\varinjlim} G_k \cong
C({\frak K},\Bbb Z)^{\infty}.
$$
\end{pf}
Therefore  we conclude
\begin{thm}
$$
K_0(\OLF) \cong {\Bbb Z} \oplus  C({\frak K},{\Bbb Z})^{\infty},
\qquad
K_1(\OLF) \cong 0.
$$
\end{thm}
\begin{pf}
Since
$
K_0(\OLF)
$
is isomorphic to
$$
{\Bbb Z}
\oplus
 \underset{l}{\varinjlim} 
\{ J_{l+1,l}:
{\Bbb Z}^{m(l)-1}
\rightarrow
{\Bbb Z}^{m(l+1)-1} \}
$$
and the second summand above denoted by $G$ is isomorphic to
$ C({\frak K},{\Bbb Z})^{\infty}$,
one gets
$
K_0(\OLF) \cong {\Bbb Z} \oplus  C({\frak K},{\Bbb Z})^{\infty}.
$
We have already seen the formula
$K_1(\OLF) \cong 0.
$
\end{pf}
Therefore Theorem 1.1 holds.


\end{document}